\documentclass[letterpaper, 10 pt, conference]{ieeeconf}

\IEEEoverridecommandlockouts                              

\usepackage{hyperref}
\usepackage{url}
\usepackage{bbold}
\usepackage{amsfonts}
\usepackage{amsmath,amssymb}
\usepackage{bbm}

\usepackage{graphicx}
\usepackage{hyperref}
\usepackage{wrapfig}
\usepackage{mathrsfs} 
\usepackage{algorithm,algorithmic}
\usepackage{array}
\usepackage{times}
\usepackage{url}
\usepackage{cite}
\usepackage{upgreek}

\usepackage{float}
\usepackage{color}

\usepackage{pgfplots}
\usepackage{color}
\usepackage{flushend}
\usepackage{multirow}
\usepackage{rotating}
\usepackage{tikz,pgfplots}
\usepackage{standalone}

\newcommand\numberthis{\addtocounter{equation}{1}\tag{\theequation}}

\newtheorem{theorem}{Theorem}

\newtheorem{assumption}{Assumption}

\newtheorem{lemma}[theorem]{Lemma}

\input{mysymbol.sty}

\begin{document}

%\title{\LARGE \bf Distributed Online Optimization in Dynamic Environments Using Mirror Descent}

%\title{\LARGE \bf On Distributed Tracking of Dynamic Parameters: Nonlinear Observation and Adversarial Noise}

\title{\LARGE \bf An Online Optimization Approach for Multi-Agent Tracking of Dynamic Parameters in the Presence of Adversarial Noise}

%\title{Mirror Descent for Decentralized Online Optimization in Dynamic Environments}

\author{Shahin Shahrampour and Ali Jadbabaie 
\thanks{This work was supported by ONR BRC Program on Decentralized,
Online Optimization.}
\thanks{Shahin Shahrampour is with the Department of Electrical Engineering at Harvard University, Cambridge, MA 02138 USA. (e-mail: {\tt shahin@seas.harvard.edu}).}
\thanks{Ali Jadbabaie is with the Institute for Data, Systems, and Society at Massachusetts Institute of Technology, Cambridge, MA 02139 USA. (email: {\tt jadbabai@mit.edu}).}
}

\maketitle

% REQUIRED      %%%%%%%%%%%%%    ABSTRACT  %%%%%%%%%%%%%%%%%
\begin{abstract}
This paper addresses tracking of a moving target in a multi-agent network. The target follows a linear dynamics corrupted by an adversarial noise, i.e., the noise is not generated from a statistical distribution. The location of the target at each time induces a global time-varying loss function, and the global loss is a sum of local losses, each of which is associated to one agent. Agents noisy observations could be nonlinear. We formulate this problem as a distributed online optimization where agents communicate with each other to track the minimizer of the global loss.
We then propose a decentralized version of the Mirror Descent algorithm and provide the
non-asymptotic analysis of the problem. Using the notion of dynamic regret, we measure the performance of our algorithm versus its offline counterpart in the centralized setting. We prove that the bound on dynamic regret scales inversely in the network spectral gap, and it represents the adversarial noise causing deviation with respect to the linear dynamics. Our result subsumes a number of results in the distributed optimization literature. Finally, in a numerical experiment, we verify that our algorithm can be simply implemented for multi-agent tracking with nonlinear observations. 
\end{abstract}

%%%%%%%%%%%%%%%%%%%%%%%%%%%%%%%%%%%%%%%%%%%%%%%

%%%%%%%%%%%%%%%%%%%%%%%    SECTIONS   %%%%%%%%%%%%%%%%

% !TEX root =  ACC-2017.tex

\section{Introduction}
Distributed estimation, detection, and tracking is ubiquitous in engineering applications ranging from sensor and robotic networks to social networks, and it has received a lot of attention for many years \cite{bullo2009distributed,kar2012distributed,shahrampour2013online,shahinjournal,nedic2016network}. In
these scenarios, the task is to estimate the value of a parameter which may or may not be dynamic. A group of agents aim to accomplish this task as a team. Each individual agent only partially
observes the parameter, but the global spread of observations in the network allows agents to estimate the parameter collaboratively. This would require agents to aggregate local information, and many methods use consensus protocols as a critical component \cite{jadbabaie2003coordination}. It is well-known that when agents' observations are linear with respect to the parameter, the tracking problem is equivalent to minimizing a {\it global} quadratic loss, written as a sum of {\it local} quadratic losses (see e.g. \cite{shahinACC2016}). However, in general, the global loss can be more complicated, resulting in  nonlinear observations.

In real-world applications, the parameter of interest is often time-varying. Therefore, regardless of the structure of the loss, the dynamic nature of the problem brings forward two issues: {\it (i)} The local losses are observed in an {\it online} or {\it sequential} fashion, i.e., the local losses are disclosed to agents only after they form their estimates at each round, and they are not aware of future loss functions. Therefore, the problem must be solved in an online setting. 
{\it (ii)} The online algorithm should mimic the performance of its offline counterpart in which the losses are known a priori. The gap between the two is often called {\it regret}. Tracking the minimizer of the global loss over time introduces the notion of {\it dynamic} regret \cite{zinkevich2003online}. This framework has been studied in centralized online optimization \cite{zinkevich2003online,hall2015online,besbes2015non,jadbabaie2015online,mokhtari2016online}, where the hardness of the problem is captured via the variation in the minimizer sequence.

To address these issues in this paper, we adopt an online optimization approach to formulate distributed tracking. We consider tracking of a dynamic parameter or a moving target in a {\it network} of agents. The dynamics of the target is linear and known to agents, but the target deviates from this dynamics due to an {\it unstructured} or {\it adversarial} disturbance or noise. In other words, the noise is not necessarily generated from a statistical distribution, or it can be highly correlated to its past values over time. At each time instance, the target induces a {\it global} convex loss whose minimizer coincides with the target location. The global loss is a sum of {\it local} losses, where each local loss is associated to a specific agent. Agents exchange noisy local gradients according to a communication protocol to track the moving target.

Our problem setup is reminiscent of a distributed Kalman \cite{olfati2007distributed}. However, we differentiate the two as follows: {\it (i)} We do not assume that the target is driven by a Gaussian noise. Nor do we assume that this noise has a statistical distribution. Instead, we consider an adversarial-noise model with {\it unknown} structure. {\it (ii)} Agents observations are not necessarily linear; in fact, the observations are noisy local gradients that are non-linear  when the loss is not quadratic. Furthermore, our focus is on the {\it finite-time} analysis rather than asymptotic results. %Our setup also differs from distributed particle filtering \cite{gu2007distributed} as it is online, and agents receive only
%one observation per iteration.

We propose a decentralized version of the Mirror Descent algorithm, developed by Nemirovksi and Yudin \cite{yudin1983problem}. Using the notion of Bregman divergence in lieu of Euclidean distance for projection, Mirror Descent has been shown to be a powerful tool in large-scale optimization. Our algorithm consists of three interleaved updates: {\it (i)} each agent follows the noisy local gradient while staying close to previous estimates in the local neighborhood; {\it (ii)} agents take into account the dynamics of the moving target; {\it (iii)} agents average their estimates in their local neighborhood in a consensus step.

We then use a {\it dynamic} notion of regret to measure the difference between our online decentralized algorithm and its offline centralized version. We establish a regret bound that scales inversely in the spectral gap of the network, and it represents the adversarial noise causing deviation with respect to the linear dynamics. We further show that from optimization perspective our result subsumes two important classes of decentralized optimization in the literature: {\it (i)} decentralized optimization of {\it time-invariant} losses, and {\it (ii)} decentralized optimization of time-variant losses for {\it fixed} targets. This generalization is achieved by allowing the loss function and the target value to vary simultaneously.  We also provide a numerical experiment to show that our algorithm can be simply implemented to work with nonlinear observations in multi-agent tracking.

\noindent
{\bf Related Literature on Decentralized Optimization:} In \cite{li2016distributed}, decentralized mirror descent has been developed for time-invariant functions in the case that agents receive the gradients with a delay. Moreover, Rabbat in \cite{rabbat2015multi} proposes a decentralized mirror descent for stochastic composite optimization problems and provide guarantees for strongly convex regularizers. Duchi et al. \cite{duchi2012dual} study dual averaging for distributed optimization, and the extension of dual averaging to online distributed optimization is considered in \cite{hosseini2013online}. Mateos-N{\'u}nez and Cort{\'e}s \cite{mateos2014distributed} consider online optimization using subgradient descent of local functions, where the graph structure is time-varying. In \cite{nedic2015decentralized}, a decentralized variant of Nesterov's primal-dual algorithm is proposed for online optimization. In \cite{7353155}, distributed online optimization is studied for strongly convex objective functions over time-varying networks. Our setup follows the work of \cite{journal} on decentralized online mirror descent, but we extend the results to {\it high probability} bounds on the dynamic regret.

% \noindent
%{\bf Related Literature on Centralized, Online Optimization in Dynamic Environments:} In dynamic online optimization, the benchmark is defined in terms of a time-varying sequence, a particular case of which is the minimizer sequence of a time-varying cost function. In \cite{zinkevich2003online}, Zinkevich develops the celebrated online gradient descent and considers its extension to time-varying sequences. The authors of \cite{hall2015online} generalize this idea to study time-varying sequences following given dynamics. Besbes et al. \cite{besbes2015non} restrict their attention to minima sequence and introduce a complexity measure for the problem in terms of variation in cost functions. For the same problem, the authors of \cite{jadbabaie2015online} develop an adaptive algorithm whose regret bound is expressed in terms of the variation of both functions and minima sequence, while in \cite{mokhtari2016online} an improved rate is derived for strongly convex objectives. %Moreover, online dynamic optimization with linear objectives is carefully discussed in \cite{lee2015resisting}. Finally, Yang et al. \cite{yang2016tracking} provide optimal bounds for when the minimizer belongs to the feasible set. 

% !TEX root =  ACC-2017.tex

\begin{center}
  \begin{tabular}{| c || l | }
    \hline
     $\vphantom{\sum^N} [n]$ &  The set $\{1,2,...,n\}$ for any integer $n$ \\ \hline 
     $\vphantom{\sum^N} x^\top$ & Transpose of the vector $x$ \\ \hline
     $\vphantom{\sum^N} x(k)$ & The $k$-th element of vector $x$ \\ \hline
     $\vphantom{\sum^N} I_n$ &  Identity matrix of size $n$ \\ \hline
     $\vphantom{\sum^N} \Delta_d$ &  The $d$-dimensional probability simplex \\ \hline
     $\vphantom{\sum^N} \inn{\cdot, \cdot}$ &  Standard inner product operator \\ \hline
     %$\vphantom{\sum^N} \e{\cdot}$ &  Expectation operator \\ \hline
     $\vphantom{\sum^N} \norm{\cdot}_p$ &  $p$-norm operator \\ \hline
     $\vphantom{\sum^N} \norm{\cdot}_*$ &  The dual norm of $\norm{\cdot}$ \\ \hline
     %$\vphantom{\sum^N} \mathbb{1}$ &  Vector of all ones \\ \hline
     $\vphantom{\sum^N} \sigma_i(W)$ &  The $i$-th largest eigenvalue of $W$ in magnitude\\ \hline
  \end{tabular}
\end{center}

\section{Problem Formulation and Algorithm}\label{prelim}
\subsection{Dynamical Model and Optimization Perspective}
Consider a $d$-dimensional moving target $x^\star_{t}$ following the linear dynamics $A$ for a {\it finite} time $T$ as
\begin{align}\label{dynamics}
x^\star_{t+1}=Ax^\star_{t}+v_t, \ \ \ \ \ \ \ \ \ \  t \in [T]
\end{align}
where $A\in \real^{d\times d}$ is {\it known}, and $v_t \in \real^d$ is an {\it adversarial} noise, i.e., the sequence $\{v_t\}_{t=1}^T$ is neither generated according to a statistical distribution, nor it is independent over time. Our goal is to track $x^\star_{t}$, and regardless of the observation model, a distribution-dependent mechanism, such as Kalman or particle filter, cannot solve the problem since the noise does not assume a statistical distribution.

In the centralized version of the tracking problem above, the observations of $x^\star_{t}$ are realized through a {\it time-varying}, {\it global} loss function. That is, consider the tracking problem above as an optimization, where $x^\star_{t}$ is the minimizer of the global loss at time $t$. Let $\X$ be a convex,  compact set, and represent the global loss by $f_t: \X \rightarrow \real$ at time $t$. As the global loss varies over time, the goal is to track the minimizer of $f_t(\cdot)$, which is $x^\star_t$. The {\it offline} and {\it centralized} version of our problem can be viewed as follows 
\begin{equation}\label{problem}
\begin{aligned}
& \underset{x_1,\ldots,x_T}{\text{minimize}}
& & \sum_{t=1}^T f_t(x_t) \\
& \text{subject to}
& & x_t \in \X, \; t \in [T].
\end{aligned}
\end{equation}
We are interested to solve the problem above in an {\it online} and {\it decentralized} fashion. In particular, the global function at time $t$ is a sum of $n$ {\it local} functions as
\begin{align}\label{decomposition}
f_t(x):=\frac{1}{n}\sum_{i=1}^nf_{i,t}(x),
\end{align}
 where $f_{i,t}: \X \rightarrow \real$ is a local convex function on $\X$ for all $i \in [n]$. We consider a network of $n$ agents facing two challenges when solving problem \eqref{problem}: {\it (i)} agent $j\in [n]$ receives information only about $f_{j,t}(\cdot)$ and does not observe the global loss function $f_t(\cdot)$, which is common to decentralized schemes; {\it (ii)} The functions are revealed to agents sequentially along the time horizon, i.e., at any time instance $s$, agent $j$ has observed $f_{j,t}(\cdot)$ for $t<s$, whereas the agent does not know $f_{j,t}(\cdot)$ for $s \leq t\leq T$, which is common to online settings. 

The agents interact with each another, and their relationship is captured via an undirected graph $\mathcal{G}=(\mathcal{V}, \mathcal{E})$, where $\mathcal{V}=[n]$ denotes the set of nodes, and $\mathcal{E}$ is the set of edges. Each agent $i$ assigns a positive weight $[W]_{ij}$ for the information received from agent $j \neq i$, and the set of neighbors of agent $i$ is defined as $\mathcal{N}_i:=\{j : [W]_{ij}>0\}$.

While the problem framework is reminiscent of a distributed Kalman \cite{olfati2007distributed}, there are fundamental distinctions in our setup: {\it (i)} The adversarial noise $v_t$ is neither Gaussian nor of known statistical distribution. It can be thought as a noise with {\it unknown} structure, which represents the deviation from the dynamics\footnote{In online optimization, the focus is not on distribution of data. Instead, data is thought to be generated arbitrarily, and its effect is observed through the loss functions\cite{shalev2011online}.}. {\it (ii)} Agents observations are not necessarily linear; in fact, the observations are local gradients of $\{f_{i,t}(\cdot)\}_{t=1}^T$ and are non-linear when the objective is not quadratic. The other implicit distinction in this work is our focus on {\it finite-time} analysis rather than asymptotic results. 

%We note that our framework also differs from distributed particle filtering\cite{gu2007distributed} since agents receive only one observation per iteration, and the adversarial noise $v_t$ has no structure or distribution. 

From optimization perspective, our framework subsumes two important classes of decentralized optimization in the literature:
\begin{itemize}
\item[1)] Existing methods often consider {\it time-invariant} objectives (see e.g. \cite{nedic2009distributed,duchi2012dual,li2016distributed}). This is simply the special case where $f_t(x)=f(x)$ and $x_t=x$ in \eqref{problem}.
\item[2)] Online algorithms deal with {\it time-varying} functions, but often the network's objective is to minimize the temporal average of $\{f_t(x)\}_{t=1}^T$ over a fixed variable $x$ (see e.g. \cite{hosseini2013online,mateos2014distributed}). This can be captured by our setup when $x_t=x$ in \eqref{problem}.
\end{itemize}
However, in the tracking problem, {\it functions} and {\it comparator variables} evolve simultaneously, i.e., the variables $\{x_t\}_{t=1}^T$ are not constrained to be fixed in \eqref{problem}. Recall that $x^\star_t:=\argmin_{x \in \X}f_t(x)$ is the minimizer of the global loss function at time $t$. Then, the solution to problem \eqref{problem} is simply $\sum_{t=1}^T f_t(x^\star_t)$. Denote by $\x_{i,t}$ the estimate of agent $i$ for $x_t^\star$ at time $t$. To exhibit the online nature of problem \eqref{problem}, we reformulate it using the notion of {\it dynamic} regret as follows
\begin{align}\label{regret}
\textbf{\textit{Reg}}^d_T=\frac{1}{n}\sum_{i=1}^n \sum_{t=1}^T f_t(\x_{i,t}) -  \sum_{t=1}^T f_t(x^\star_t).
\end{align}
Then, the objective is to minimize the dynamic regret above which measures the gap between the online algorithm and its offline version. Our performance bound shall exhibit the impact of system noise, i.e., we want to prove a regret bound in terms of 
\begin{align}\label{CT}
\norm{v_t}=\norm{x^\star_{t+1}-Ax^\star_{t}},
\end{align}
which represents the deviation of the moving target with respect to dynamics $A$. Note that generalizing the results to the linear time-variant dynamics is straightforward, i.e., when $A$ is replaced by $A_t$ in \eqref{dynamics}.

 \subsection{Technical Assumptions}

To solve the multi-agent online optimization \eqref{regret}, we propose to decentralize the Mirror Descent algorithm \cite{yudin1983problem}. Mirror Descent has been shown to be a powerful method in large-scale optimization by using Bregman divergence in lieu of Euclidean distance in the projection step. Before defining Bregman divergence and elaborating the algorithm, we start by stating a couple of standard assumptions on loss functions and agents communication.

\begin{assumption}\label{A1}
For any $i\in [n]$, the function $f_{i,t}(\cdot)$ is Lipschitz continuous on $\X$ with a uniform constant $L$. That is, $$|f_{i,t}(x)-f_{i,t}(y)|\leq L\norm{x-y},$$ for any $x,y \in \X$. %This further implies that the gradient of $f_{i,t}(\cdot)$ denoted by $\nabla f_{i,t}(\cdot)$ is uniformly bounded on $\X$ by the constant $L$, i.e., we have $\norm{\nabla f_{i,t}(\cdot)}_* \leq L$.\footnote{This relationship is standard, see e.g. Lemma 2.6. in \cite{shalev2011online} for more details.}
\end{assumption}

\begin{assumption}\label{A2}
The network is connected\footnote{The setup is generalizable to when network connectivity changes over time, and the communication matrix is time-varying.}, i.e., there exists a path from any agent $i\in [n]$ to any agent $j\in [n]$. Also, the matrix $W$ is symmetric and doubly stochastic with positive diagonal. That is,
$$
\sum_{i=1}^n[W]_{ij}=\sum_{j=1}^n[W]_{ij}=1.
$$
\end{assumption}
The connectivity constraint in Assumption \ref{A2} guarantees the information flow in the
network. 

We now outline the notion of Bregman divergence, which is critical in the development of Mirror Descent. Consider a compact, convex set $\X$, and let $\R : \X \rightarrow \real$ denote a 1-strongly convex function on $\X$ with respect to a norm $\norm{\cdot}$. That is,
$$
\R(x)\geq \R(y)+\inn{\nabla\R(y),x-y}+\frac{1}{2}\norm{x-y}^2.
$$
for any $x,y \in \X$. Then, the Bregman divergence $\dr(\cdot,\cdot)$ with respect to the function $\R(\cdot)$ is defined as follows:
$$
\dr(x,y):=\R(x)-\R(y)-\inn{x-y,\nabla\R(y)}.
$$ 
The definition of the Bregman divergence and the strong convexity of $\R(\cdot)$ imply that
\begin{align}\label{bregcond}
\dr(x,y) \geq \frac{1}{2}\norm{x-y}^2,
\end{align}
for any $x,y \in \X$. Two famous examples of Bregman divergence are the Euclidean distance and the Kullback-Leibler (KL) divergence generated from $\R(x)=\frac{1}{2}\norm{x}^2_2$ and $\R(x)=\sum_{i=1}^dx(i)\log x(i)-x(i)$, respectively.

\begin{assumption}\label{A3}
Let $x$ and $\{y_i\}_{i=1}^n$ be vectors in $\real^d$. We assume that the Bregman divergence satisfies the separate convexity in the following sense 
$$\dr(x,\sum_{i=1}^n\alpha(i)y_i) \leq \sum_{i=1}^n\alpha(i)\dr(x,y_i),$$
where $\alpha\in\Delta_n$ is on the $n$-dimensional simplex.
\end{assumption}
The assumption is satisfied for commonly used cases of Bregman divergence. For instance, the Euclidean distance evidently respects the condition. The KL-divergence also satisfies the constraint,  and we refer the reader to Theorem 6.4. in \cite{bauschke2001joint} for the proof.
\begin{assumption}\label{A4}
The Bregman divergence satisfies a Lipschitz condition of the form
\begin{align*}
| \dr(x,z)-\dr(y,z) | \leq K\|x- y\|, 
\end{align*}
for all $x,y,z \in \X$.
\end{assumption}
When the function $\R$ is Lipschitz on $\X$, the Lipschitz condition on the Bregman divergence is automatically satisfied. Again, for the Euclidean distance the assumption evidently holds. In the particular case of KL divergence, the condition can be achieved via mixing a uniform distribution to avoid the boundary
 (see e.g. \cite{jadbabaie2015online} for more comments on the assumption).

\begin{assumption}\label{A5}
The dynamics $A$ is assumed to be {\it non-expansive}. That is, the condition
\begin{align*}
\dr\big(Ax,Ay\big) \leq \dr\big(x,y\big),
\end{align*}
holds for all $x,y \in \X$, and $\norm{A} \leq 1$.
\end{assumption} 
The assumption postulates a natural constraint on the dynamics $A$: it does not allow the effect of a poor estimation (at one step) to be amplified as the algorithm moves forward.

 \subsection{Decentralized Tracking via Online Mirror Descent}
We now propose our algorithm to solve the problem formulated in terms of dynamic regret in \eqref{regret}. In our setting, agents observations are gradients of the local losses. However, common in distributed state estimation and tracking, these observations are noisy. Hence, denoting the local gradient of agent $i$ at time $t$ by $\nabla_{i,t}:=\nabla f_{i,t}(\x_{i,t})$, the agent only receives ${\boldsymbol \nabla}_{i,t}$ representing the stochastic gradient. The stochastic oracle that provides noisy gradients satisfies the following constraints\footnote{For simplicity, we use one constant $L$ to bound gradients as well as the stochastic gradients.} 
\begin{align}\label{condition}
\e{\vphantom{\norm{{\boldsymbol \nabla}_{i,t}}_*^2}{\boldsymbol \nabla}_{i,t}\big\vert \F_{t-1}}=\nabla_{i,t} \ \ \ \ \ \ \ \ \norm{{\boldsymbol \nabla}_{i,t}}_* \leq L,
\end{align}
where $\F_t$ is the $\sigma$-field containing all information prior to the outset of round $t+1$.
A commonly used model to generate stochastic gradients satisfying \eqref{condition} is an additive, bounded, zero-mean noise. Agents then track the moving target using a decentralized variant of Mirror Descent as follows\footnote{We set $\x_{i,t}$ to be the vector of all zeros to initialize the algorithm. In
general, any initialization could work for the algorithm.}
\begin{subequations}
\begin{align}
\xhx_{i,t+1}&=\argmin_{x\in \X}  \big\{ \eta_t\inn{x, {\boldsymbol \nabla}_{i,t}} + \dr(x,\y_{i,t}) \big\}, \label{xhupdate2}\\
\x_{i,t}&=A \xhx_{i,t}, \ \ \ \  \text{and}  \ \ \ \ \  \y_{i,t}=\sum_{j=1}^n [W]_{ij} \x_{j,t} \label{xyupdate2},
\end{align}
\end{subequations}
where $\{\eta_t\}_{t=1}^T$ is the step-size sequence, and $A\in \real^{d \times d}$ is the given dynamics in \eqref{dynamics} which is {\it common} knowledge. In these updates, $\x_{i,t}\in \real^d$ represents the estimate of agent $i$ of the moving target $x^\star_t$ at time $t$. The step-size sequence should be tuned for different cases, but it is generally non-increasing and positive. 

The update \eqref{xhupdate2} allows an agent to follow the noisy local gradient while keeping the estimate close to those of the local neighborhood. This closeness occurs by minimizing the Bregman divergence. On the other hand, the first update in \eqref{xyupdate2} takes into account the dynamics of the moving target, and the second update in \eqref{xyupdate2} is the {\it consensus} term averaging the estimates in the local neighborhood.

% !TEX root =  ACC-2017.tex

\section{Theoretical Results}\label{theory}
In this section, we state our theoretical result on the {\it non-asymptotic} performance of the decentralized online mirror descent for tracking dynamic parameters. Theorem \ref{theorem1} proves a bound on the dynamic regret, which captures the deviation of the moving target from the dynamics $A$ (tracking error), the decentralization cost (network error), and the impact of stochastic gradients (stochastic error). We show that this theorem recovers previous rates on decentralized optimization once the tracking error is removed. Also, it recovers previous rates on centralized online optimization in dynamic setting when the network error is eliminated. The proof is given in Appendix (Section \ref{appendix}).
\begin{theorem}\label{theorem1}{\it
Consider a moving target $x_t^\star \in \real^d$ with the dynamical model of \eqref{dynamics}. Further consider the distributed, online tracking problem formulated in \eqref{regret}, where $\x_{i,t}$ denotes the local estimate of agent $i\in [n]$ of the moving target $x^\star_t$ at time $t\in [T]$. Let the local estimates be generated by updates \eqref{xhupdate2}-\eqref{xyupdate2}, where the stochastic gradients satisfy the condition \eqref{condition}. Given Assumptions [\ref{A1}-\ref{A5}], the dynamic regret can be bounded as}
\begin{align*}
\textbf{\textit{Reg}}^d_T &\leq {\tt E_{Track}}+{\tt E_{Net}}+{\tt E_{Stoch}},
\end{align*}
{\it with probability at least $1-\delta$, where}
\begin{align*}
\vphantom{\sum_{t=1}^T\frac{K}{\eta_{t+1}}}{\tt E_{Track}}&:=\frac{2R^2}{\eta_{T+1}}+\sum_{t=1}^T\frac{K}{\eta_{t+1}}\norm{x^\star_{t+1}-Ax^\star_{t}}+L^2\sum_{t=1}^T\frac{\eta_t}{2}\\
\vphantom{\sum_{t=1}^T\frac{K}{\eta_{t+1}}}{\tt E_{Net}}&:=4L^2\sqrt{n}\sum_{t=1}^T\sum_{\tau=0}^{t-1}\eta_\tau\sigma^{t-\tau-1}_2(W)\\
\vphantom{\sum_{t=1}^T\frac{K}{\eta_{t+1}}}{\tt E_{Stoch}}&:= 8LR \sqrt{-T\log \delta},
\end{align*}
{\it and $R^2:=\sup_{x,y\in \X} \dr(x, y).$}
\end{theorem}
In view of \eqref{dynamics}, the dynamical model of the target is described with the noise $v_t$. The term ${\tt E_{Track}}$ shows the dependence of performance bound to noise by aggregating the errors $\norm{x^\star_{t+1}-Ax^\star_{t}}=\norm{v_t}$ over time. Also, ${\tt E_{Net}}$ and ${\tt E_{Stoch}}$ are the errors related to network and stochastic gradients, respectively. 

In Section \ref{prelim}, we discussed that our setup generalizes some of the previous results. It is now important to see that this generalization is valid in the sense that our result can recover those special cases:
\begin{itemize}
\item[$\triangleright$] When the global loss $f_t(x)=f(x)$ is time-invariant, the target $\{x^\star_t\}_{t=1}^T$ is fixed, i.e., the dynamics $A=I_d$ and $v_t=\mathbb{0}$ in \eqref{dynamics}. In this case in Theorem \ref{theorem1}, the term involving $\norm{x^\star_{t+1}-Ax^\star_{t}}$ in ${\tt E_{Track}}$ is equal to zero, and we can use the step-size sequence $\eta=\sqrt{(1-\sigma_2(W))/T}$ to recover the result of comparable algorithms, such as Theorem 4 in \cite{duchi2012dual} on distributed dual averaging. 
\item[$\triangleright$] The same argument holds when the global loss is time-variant, but the target is fixed. This setup is studied, for instance, in \cite{hosseini2013online} via distributed online dual averaging with exact gradients. Disregarding ${\tt E_{Stoch}}$ in our bound due to stochastic gradients, since $\norm{x^\star_{t+1}-Ax^\star_{t}}=0$ again, we recover Corollary 3 in \cite{hosseini2013online}.  
\item[$\triangleright$] When the graph is complete, $\sigma_2(W)=0$ and hence ${\tt E_{Net}}=0$. We then recover the results of  \cite{hall2015online} on centralized online learning (for linear dynamics) with exact gradients once we remove ${\tt E_{Stoch}}$ due to stochastic gradients. 
\end{itemize}
%When mismatch errors $\{v_t\}_{t=1}^T$ are large, the minimizer sequence $\{x^\star_t\}_{t=1}^T$ oscillates wildly, and $C_T$ can become linear in time. The bound in the corollary is then not useful for keeping the dynamic regret sub-linear. Such behavior is natural since even in the centralized online optimization, the algorithm receives only a {\it single} gradient to predict the next step. As discussed in Section \ref{prelim}, in this worst-case, the problem is generally intractable. Our primary goal, however, was to consider $C_T$ as a complexity measure of the problem environment and express the regret bound with respect to this parameter. In practice, if the algorithm is allowed to query {\it multiple} gradients per time, the error would be reduced, but this direction is beyond the scope of this paper. 

% !TEX root =  ACC-2017.tex

\section{Numerical Experiment: Tracking Maneuvering Targets}\label{application}
In Mirror Descent algorithm, one has freedom over the selection of the Bregman divergence. A particularly well-known type of Bregman is the Euclidean distance, commonly used in state estimation and tracking dynamic parameters. We focus on this scenario in this section to provide the numerical experiments for our method. 

We consider a slowly maneuvering target in the $2D$ plane and assume that each position component of the target evolves independently according to a near constant velocity model \cite{bar1987tracking}. The state of the target at each time consists of four components: horizontal position, vertical position, horizontal velocity, and vertical velocity. We represent the state at time $t$ by $x^\star_{t} \in \real^{4}$, and therefore, the state space model takes the form
$$
x^\star_{t+1}=A x^\star_{t} + v_t,
$$ 
where $v_t \in \real^4$ is the system noise, and using $\otimes$ for Kronecker product, $A$ can be written as
$$
A=I_2 \otimes \begin{bmatrix}
1 & \epsilon \\
0 & 1
\end{bmatrix},
$$
with $\epsilon$ being the sampling interval\footnote{The sampling interval of $\epsilon$ (seconds) is equivalent to the sampling rate of $1/\epsilon~(Hz)$.}. The goal is to cooperatively track $x^\star_{t}$ in a network of agents. This problem has been studied in the context of distributed Kalman filtering \cite{olfati2007distributed,cattivelli2010diffusion}, state estimation \cite{khan2010connectivity,das2013distributed,han2015stochastic}, and particle filtering \cite{gu2007distributed,hlinka2012likelihood,li2015distributed}. However, in contrast to Kalman filtering, we do not assume that the system noise $v_t$ is Gaussian. Also, as opposed to particle filtering, we do not receive a large number of samples (particles) per iteration since our setup is online, i.e., agents only observe one sample per time. Furthermore, we do not assume a statistical distribution on $v_t$ in our analysis, which differentiates our framework from state estimation. We adopt a model-free approach where the noise can be adversarial (deterministic), stochastic with dependence over time, or of some complex structure. We generate the noise as follows. At each time $t$ we draw a sample $\nu_t \in \real^{4}$ from a zero-mean Gaussian distribution with covariance matrix $\Sigma$ as follows
$$
\Sigma=\sigma^2_{\nu} I_2 \otimes \begin{bmatrix}
\epsilon^3/3 & \epsilon^2/2 \\
\epsilon^2/2 & \epsilon
\end{bmatrix},
$$ 
for the sampling interval $\epsilon=0.1$ seconds which amounts to frequency $10~Hz$. Then, we let the system noise be $v_t=\nu_t\norm{\nu_t}_\infty$. Though $\nu_t$ is generated from Gaussian distribution, the mismatch noise $v_t$ is non-Gaussian and can have a complicated distribution. The constant $\sigma^2_\nu$ takes different values in each experiment, and we describe this choice later.

We consider a sensor network of $n=25$ agents located on a $5 \times 5$ grid. Agents aim to track the moving target $x^\star_{t}$ collaboratively. Agents observe a noisy version of the target through a {\it local} loss function, and these observations are nonlinear. In particular, let the quantity $\z_{i,t}$ be a noisy version of one coordinate of $x^\star_{t}$ as follows
$$
\z_{i,t}=\ee_{k_i}^\top x^\star_{t} + \w_{i,t},
$$
where $\w_{i,t} \in \real$ denotes a random noise, and $\ee_k$ is the $k$-th unit vector in the standard basis of $\real^4$ for $k\in \{1,2,3,4\}$. We partition the agents into four groups, and for each group we select one specific $k_i$ from the set $\{1,2,3,4\}$. The random noise $\w_{i,t}$ satisfies the standard assumption of being zero-mean and finite-variance. Again, to show that our results are not dependent on Gaussian noise, we generate $\w_{i,t}$ independently from a uniform distribution on $[-1,1]$.

Then, at time $t$ the {\it local} loss for agent $i$ takes the form
$$
f_{i,t}(x):=\frac{1}{4}\e{\left(\z_{i,t}-\ee_{k_i}^\top x\right)^4  \big\vert \F_{t-1},x^\star_{t}},
$$
resulting in the {\it global} loss  
\begin{align*}
f_t(x):=\frac{1}{4n}\sum_{i=1}^n\e{\left(\z_{i,t}-\ee_{k_i}^\top x\right)^4 \big\vert \F_{t-1},x^\star_{t}},
\end{align*}
where $\F_{t-1}$ is the $\sigma$-field containing all information in $\{\w_{i,s}\}_{s=1}^{t-1}$. It is straightforward to see that $x^\star_{t}$ is the minimizer of the global loss. Observation of agent $i$ at time $t$ is the stochastic gradient of the local loss
$$
\nabla f_{i,t}(x)=\left(\z_{i,t}-\ee_{k_i}^\top x\right)^3 \ee_{k_i}.
$$
We derive an explicit update to form an estimate $\x_{i,t}$ of $x^\star_{t}$. We use Euclidean distance as the Bregman divergence in updates \eqref{xhupdate2}-\eqref{xyupdate2} to get\footnote{We assume that the state of the target remains in a convex, compact set, and the updates can keep the estimate in the set without the projection step. This assumption can be satisfied in the finite-time domain.} 
\begin{align*}
\x_{i,t}&=\sum_{j=1}^n [W]_{ij} A\x_{j,t-1}+ \eta_t A\ee_{k_i} \left(\z_{i,t-1}-\ee_{k_i}^\top \x_{i,t-1}\right)^3,
\end{align*}
and tune the step size to $\eta_t=\eta=0.1$. % since using diminishing step size is not useful in tracking unless we have diminishing system noise \cite{acemoglu2008convergence}. 
 The update is akin to {\it consensus+innovation} updates in the literature (see e.g. \cite{kar2012distributed,shahinACC2016}) though we recall that the observation is nonlinear, and the system noise $v_t$ is arbitrary.  
\begin{figure}[t!]
\centering
\includegraphics[trim = 10mm 53mm 0mm 58mm, clip, scale=0.45]{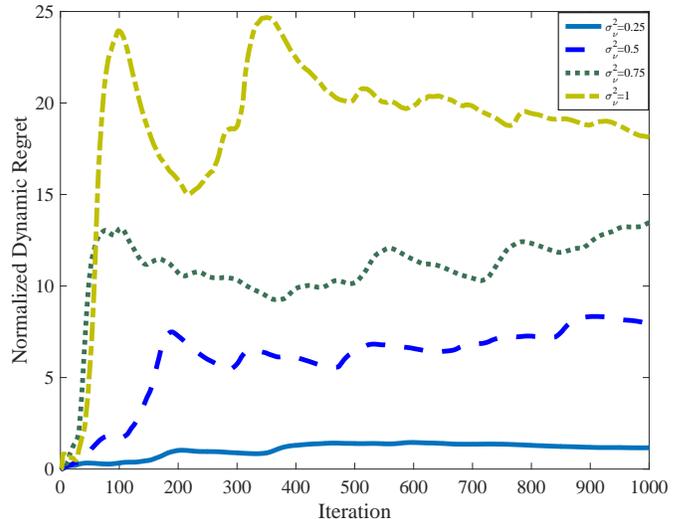}
\caption{The plot of dynamic regret versus iterations. Naturally, when $\sigma^2_v$ is smaller, the innovation noise added to the dynamics is smaller with high probability, and the network incurs a lower dynamic regret. In this plot, the dynamic regret is normalized by iterations, so the $y$-axis is $\textbf{\textit{Reg}}^d_T/T$.}
\label{Regretvs}
\end{figure}

It is proved in \cite{shahinACC2016} that in decentralized tracking, the dynamic regret can be presented in terms of the tracking error $\x_{i,t}-x^\star_{t}$ when the local losses are quadratic. More specifically, the expected dynamic regret averages the tracking error over space and time (when normalized by $T$). While here we deal with polynomial loss of power four, the connection between tracking error and dynamic regret still holds true. Therefore, using the result of Theorem \ref{theorem1} we can expect that once the parameter does not deviate too much from the dynamics, i.e., when $\sum_{t=1}^T\norm{v_t}$ is small, the bound on the dynamic regret as well as the collective tracking error is small. 

\begin{figure}[h!]
\centering
\includegraphics[trim =21mm 52mm 0mm 47mm, clip, scale=0.47]{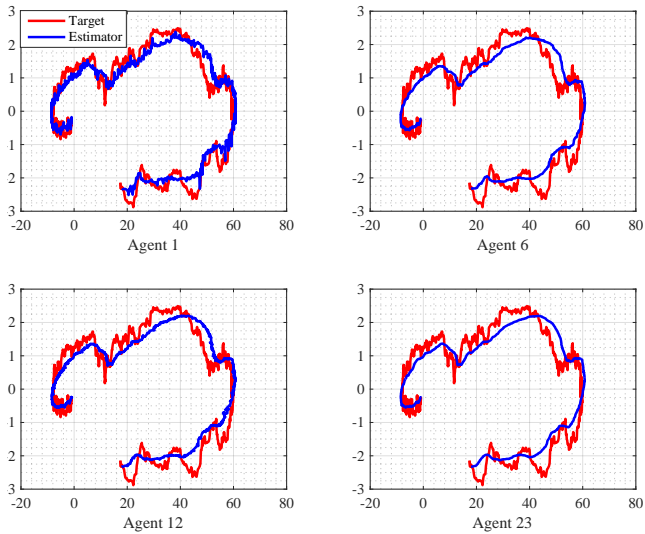}
\caption{The trajectory of $x^\star_t$ over $T=1000$ iterations is shown in red. We also depict the trajectory of the estimator $\x_{i,t}$ (shown in blue) for $i\in \{1,6,12,23\}$ and observe that it closely follows $x^\star_t$ in every case.}
\label{Agents}
\end{figure}

We show this intuitive idea by setting $\sigma^2_\nu$ to different vaues. Larger values for $\sigma^2_\nu$ are expected to cause more deviations from the dynamics $A$ and larger dynamic regret (worse performance). In Fig. \ref{Regretvs}, we plot the normalized dynamic regret for $\sigma^2_\nu\in \{0.25,0.5,0.75,1\}$. Note that for each value of $\sigma^2_\nu$, we run the experiment only once to investigate the high probability bound in Theorem \ref{theorem1}. As expected, the performance improves once $\sigma^2_v$ tends to smaller values.

We next restrict our attention to the case that $\sigma^2_v=0.5$. For one run of this case, we provide a snapshot of the target trajectory (in red) in Fig. \ref{Agents} and plot the estimator trajectory (in blue) for agents $i\in \{1,6,12,23\}$. Fig. \ref{Agents} suggests that agents' estimators closely follow the trajectory of the moving target with high probability.

%%%%%%%%%%%%%%%%%%%%%%%%%%%%%%%%%%%%%%%%%%%%%%%

%%%%%%%%%%%%%%%%%%%%%%%   CONCLUSION     %%%%%%%%%%%%%%
\section{Conclusion}\label{conclusion}
In this paper, we addressed tracking of a moving target in a network of agents. The target follows a linear dynamics which is common knowledge to agents, but it deviates from this dynamics due to an additive noise of an unknown structure. We formulated the problem as an online optimization of a global time-varying loss in a distributed fashion. The global loss at each time is a sum of a finite number of local losses, and each agent in the network holds a private copy of one local loss. Agents are unaware of the future loss functions as the local losses only become available to them sequentially. They exchange noisy local gradients with each other to track the value of the target.

Our proposed algorithm for this setup can be cast as a decentralized version of Mirror Descent. We however incorporated two more steps to include agents interactions and dynamics of the target. We used a notion of network dynamic regret to measure the performance of our algorithm versus its offline counterpart. We established that the regret bound scales inversely in the spectral gap of the
network and captures the deviation of the target with respect to the dynamics.  Our results generalized a number of results in online and offline distributed optimization. Also, numerical experiments verified the applicability of our algorithm to multi-agent tracking with nonlinear observations. Future directions include studying the algorithm in the case that several observations are available per round, i.e., when agents can receive multiple noisy gradients per time. The method can be useful in the sensor networks where each sensor can have multiple measurements from different sources.

%%%%%%%%%%%%%%%%%%%%%%%%%%%%%%%%%%%%%%%%%%%%%%%

%%%%%%%%%%%%%%%%%%%%%%%%%%    APEENDIX    %%%%%%%%%%%%%

%%%%%%%%%%%%%%%%%%%%%%%%%%%%%%%%%%%%%%%%%%%%%%%

\bibliographystyle{IEEEtran}
\bibliography{IEEEabrv,references}

% !TEX root =  ACC-2017.tex

\section{Appendix}\label{appendix}
We make use of two technical lemmas (Lemma \ref{meandeviation} and \ref{auxlemma}) proved in the Appendix of \cite{journal}. We state their results here and use them in the proof of Theorem \ref{theorem1}.

\begin{lemma}\label{meandeviation}
{\it Let $\X$ be a convex set in a Banach space $\B$, $\R : \B \rightarrow \real$ denote a 1-strongly convex function on $\X$ with respect to a norm $\|\cdot\|$, and $\dr(\cdot,\cdot)$ represent the Bregman divergence with respect to $\R$, respectively. Furthermore, assume that the local functions are Lipschitz continuous (Assumption \ref{A1}), the matrix $W$ is doubly stochastic (Assumption \ref{A2}), and the mapping $A$ is non-expansive (Assumption \ref{A5}). Then, the local estimates $\{\x_{i,t}\}_{t=1}^T$ generated by the updates \eqref{xhupdate2}-\eqref{xyupdate2} satisfy 
 $$\norm{\x_{i,t+1}-\xbx_{t+1}}\leq L\sqrt{n}\sum_{\tau=0}^t\eta_\tau\sigma^{t-\tau}_2(W),$$
 for any $i\in [n]$, where $\xbx_t:=\frac{1}{n}\sum_{i=1}^n\x_{i,t}$.}
\end{lemma}

\begin{lemma}\label{auxlemma}
{\it Let $\X$ be a convex set in a Banach space $\B$, $\R : \B \rightarrow \real$ denote a 1-strongly convex function on $\X$ with respect to a norm $\|\cdot\|$, and $\dr(\cdot,\cdot)$ represent the Bregman divergence with respect to $\R$, respectively. Furthermore, assume that the matrix $W$ is doubly stochastic (Assumption \ref{A2}), the Bregman divergence satisfies the Lipschitz condition and the separate convexity (Assumptions \ref{A3}-\ref{A4}), and the mapping $A$ is non-expansive (Assumption \ref{A5}). Then, it holds that}
\begin{align*}
\frac{1}{n}\sum_{i=1}^n\sum_{t=1}^T\bigg(\frac{1}{\eta_t}&\dr(x^\star_t,\y_{i,t})-\frac{1}{\eta_t}\dr(x^\star_t,\xhx_{i,t+1})\bigg) \\
&\leq \frac{2R^2}{\eta_{T+1}}+\sum_{t=1}^T\frac{K}{\eta_{t+1}}\norm{x^\star_{t+1}-Ax^\star_{t}},
\end{align*}
{\it where $R^2:=\sup_{x,y\in \X} \dr(x, y).$}
\end{lemma}

In what follows, we provide the proof of Theorem \ref{theorem1}. 

\subsection{Proof of Theorem \ref{theorem1}}
Recall the definition of dynamic regret in \eqref{regret}. 
Using the Lipschitz continuity of $f_t(\cdot)$ (Assumption \ref{A1}) as well as the fact that the global loss is the sum of local losses (Eq. \eqref{decomposition}), we get
\begin{align*}
\vphantom{\frac{1}{n}\sum_{i=1}^n}f_t(\x_{i,t})& - f_t(x^\star_t)=f_t(\x_{i,t})-f_t(\xbx_t)+f_t(\xbx_t) - f_t(x^\star_t)\\
\vphantom{\frac{1}{n}\sum_{i=1}^n}&\leq L\norm{\x_{i,t}-\xbx_t}+f_t(\xbx_t)-f_t(x^\star_t)\\
\vphantom{\frac{1}{n}\sum_{i=1}^n}&=L\norm{\x_{i,t}-\xbx_t}+\frac{1}{n}\sum_{i=1}^nf_{i,t}(\xbx_t)-\frac{1}{n}\sum_{i=1}^nf_{i,t}(x^\star_t),
\end{align*}
Using the Lipschitz continuity of $f_{i,t}(\cdot)$ for $i\in [n]$, we simplify above as follows  
\begin{align*}
&f_t(\x_{i,t}) - f_t(x^\star_t) \leq \frac{1}{n}\sum_{i=1}^nf_{i,t}(\x_{i,t})-\frac{1}{n}\sum_{i=1}^nf_{i,t}(x^\star_t)\\
\vphantom{\frac{1}{n}\sum_{i=1}^n}&~~~~~~~~~~~~~+L\norm{\x_{i,t}-\xbx_t}+\frac{L}{n}\sum_{i=1}^n\norm{\x_{i,t}-\xbx_t}. \numberthis \label{eq11}
\end{align*}
The second line can be controlled via Lemma \ref{meandeviation}, so we focus on the first term in the above bound. We have by convexity of $f_{i,t}(\cdot)$ that
\begin{align*}
&f_{i,t}(\x_{i,t})-f_{i,t}(x^\star_t) \leq \inn{\nabla_{i,t}, \x_{i,t}-x^\star_t}\\
&~~= \inn{{\boldsymbol \nabla}_{i,t}, \x_{i,t}-x^\star_t}+\inn{\nabla_{i,t}-{\boldsymbol \nabla}_{i,t}, \x_{i,t}-x^\star_t}\\
&~~=\inn{{\boldsymbol \nabla}_{i,t}, \xhx_{i,t+1}-x^\star_t}+\inn{{\boldsymbol \nabla}_{i,t}, \x_{i,t}-\y_{i,t}}\\
&~~+\inn{{\boldsymbol \nabla}_{i,t}, \y_{i,t}-\xhx_{i,t+1}}+\inn{\nabla_{i,t}-{\boldsymbol \nabla}_{i,t}, \x_{i,t}-x^\star_t},  \numberthis \label{eq6}
\end{align*}
We now need to bound each of the terms on the right hand side of \eqref{eq6}. The stochastic gradients are bounded in view of \eqref{condition}. Therefore, using H\"{o}lder's inequality for any primal-dual norm pair, we get 
\begin{align*}
\vphantom{\frac{1}{2\eta_t}\norm{y_{i,t}-\hx_{i,t+1}}^2+\frac{\eta_t}{2}L^2}   \inn{{\boldsymbol \nabla}_{i,t}, \y_{i,t}-\xhx_{i,t+1}} &\leq \norm{\y_{i,t}-\xhx_{i,t+1}}\norm{{\boldsymbol \nabla}_{i,t}}_*\\
\vphantom{\frac{1}{2\eta_t}\norm{y_{i,t}-\hx_{i,t+1}}^2+\frac{\eta_t}{2}L^2} &\leq L\norm{\y_{i,t}-\xhx_{i,t+1}}\\
\vphantom{\frac{1}{2\eta_t}\norm{y_{i,t}-\hx_{i,t+1}}^2+\frac{\eta_t}{2}L^2} &\leq \frac{1}{2\eta_t}\norm{\y_{i,t}-\xhx_{i,t+1}}^2+\frac{\eta_t}{2}L^2, \numberthis \label{eq7}
\end{align*}
where the last line is due to AM-GM inequality. We now recall update \eqref{xyupdate2} and use Assumptions \ref{A1} and \ref{A2} to derive 
\begin{align*}
\vphantom{2G^2\sqrt{n}\sum_{\tau=0}^t\sigma^{t-\tau}_2(W)\eta_\tau} &\inn{{\boldsymbol \nabla}_{i,t}, \x_{i,t}-\y_{i,t}}=\inn{{\boldsymbol \nabla}_{i,t}, \x_{i,t}-\xbx_{t} + \xbx_{t}-\y_{i,t}}\\
\vphantom{2G^2\sqrt{n}\sum_{\tau=0}^t\sigma^{t-\tau}_2(W)\eta_\tau} &~~~~~~~~~~~~~~=\inn{{\boldsymbol \nabla}_{i,t}, \x_{i,t}-\xbx_{t}}+\sum_{j=1}^n[W]_{ij}\inn{{\boldsymbol \nabla}_{i,t},\xbx_{t}-\x_{j,t}}\\
\vphantom{2G^2\sqrt{n}\sum_{\tau=0}^t\sigma^{t-\tau}_2(W)\eta_\tau} &~~~~~~~~~~~~~~\leq L\norm{\x_{i,t}-\xbx_{t}}+L\sum_{j=1}^n[W]_{ij}\norm{\x_{j,t}-\xbx_{t}}\\
\vphantom{2G^2\sqrt{n}\sum_{\tau=0}^t\sigma^{t-\tau}_2(W)\eta_\tau} &~~~~~~~~~~~~~~\leq 2L^2\sqrt{n}\sum_{\tau=0}^{t-1}\eta_\tau\sigma^{t-\tau-1}_2(W), \numberthis \label{eq8}
\end{align*}
where in the last line we appealed to Lemma \ref{meandeviation}. Finally, the optimality of $\xhx_{i,t+1}$ in \eqref{xhupdate2} implies (see e.g. Lemma 4.1 in \cite{beck2003mirror}) that
\begin{align*}
\inn{{\boldsymbol \nabla}_{i,t}, \xhx_{i,t+1}-x^\star_t}& \leq \frac{1}{\eta_t}\dr(x^\star_t,\y_{i,t})-\frac{1}{\eta_t}\dr(x^\star_t,\xhx_{i,t+1})\\
&-\frac{1}{\eta_t}\dr(\xhx_{i,t+1},\y_{i,t})\\
& \leq \frac{1}{\eta_t}\dr(x^\star_t,\y_{i,t})-\frac{1}{\eta_t}\dr(x^\star_t,\xhx_{i,t+1})\\
&-\frac{1}{2\eta_t}\norm{\xhx_{i,t+1}-\y_{i,t}}^2, \numberthis \label{eq9}
\end{align*}
since the Bregman divergence satisfies $\dr(x,y) \geq \frac{1}{2}\norm{x-y}^2$ for any $x,y \in \X$ in view of \eqref{bregcond}. Substituting \eqref{eq7}, \eqref{eq8}, and \eqref{eq9} into the bound \eqref{eq6}, we derive 
%\begin{align*}
%f_{i,t}(x_{i,t})-f_{i,t}(x^\star_t)&\leq \frac{\eta_t}{2}L^2+2L^2\sqrt{n}\sum_{\tau=0}^{t-1}\eta_\tau\sigma^{t-\tau-1}_2(W)\\
%&+\frac{1}{\eta_t}\dr(x^\star_t,y_{i,t})-\frac{1}{\eta_t}\dr(x^\star_t,\hx_{i,t+1}). \numberthis \label{eq14}
%\end{align*}
%Summing over $t\in [T]$ and $i\in [n]$, and applying Lemma \ref{auxlemma} completes the proof. 
%Moreover, as in Lemma \ref{meandeviation}, any bound involving $L$ which was originally an upper bound on the exact gradient must be replaced by the norm of stochastic gradient, which changes inequality \eqref{eq14} to 
\begin{align*}
f_{i,t}(\x_{i,t})-f_{i,t}(x^\star_t)&\leq \frac{\eta_t}{2}L^2+2L^2\sqrt{n}\sum_{\tau=0}^{t-1}\eta_\tau\sigma^{t-\tau-1}_2(W)\\
\vphantom{\sum_{i=1}^n}&+\frac{1}{\eta_t}\dr(x^\star_t,\y_{i,t})-\frac{1}{\eta_t}\dr(x^\star_t,\xhx_{i,t+1})\\
\vphantom{\sum_{i=1}^n}&+\inn{\nabla_{i,t}-{\boldsymbol \nabla}_{i,t}, \x_{i,t}-x^\star_t}. \numberthis \label{eq101}
\end{align*}
To bound the last term, we note that
\begin{align*}
%&\e{\inn{\nabla_{i,t}-{\boldsymbol \nabla}_{i,t}, \x_{i,t}-x^\star_t}}\\
&\e{\inn{\nabla_{i,t}-{\boldsymbol \nabla}_{i,t}, \x_{i,t}-x^\star_t} \big\vert \F_{t-1}}\\
&~~~~~~~~~~~=\inn{\e{\nabla_{i,t}-{\boldsymbol \nabla}_{i,t}\big\vert \F_{t-1}}, \x_{i,t}-x^\star_t}=0. 
\end{align*}
Also, due to \eqref{bregcond} we have $\norm{\x_{i,t}-x^\star_t}^2\leq 2\dr(\x_{i,t},x^\star_t) \leq 2R^2$, which entails {\small
$$
 \inn{\nabla_{i,t}-{\boldsymbol \nabla}_{i,t}, \x_{i,t}-x^\star_t} \leq \norm{\x_{i,t}-x^\star_t}\norm{\nabla_{i,t}-{\boldsymbol \nabla}_{i,t}}_* \leq4LR.
$$}
Therefore, summing $\inn{\nabla_{i,t}-{\boldsymbol \nabla}_{i,t}, \x_{i,t}-x^\star_t}$ over $t \in [T]$ forms a bounded difference martingale, and we can use Azuma's inequality to get
$$
\pr\left(\sum_{i=1}^n\sum_{t=1}^T\inn{\nabla_{i,t}-{\boldsymbol \nabla}_{i,t}, \x_{i,t}-x^\star_t} \geq \varepsilon\right) \leq e^{-\frac{\varepsilon^2}{32Tn^2L^2R^2}}.
$$ 
Setting the probability above to $\delta$ and solving for $\varepsilon$ implies
$$
\frac{1}{n}\sum_{i=1}^n\sum_{t=1}^T\inn{\nabla_{i,t}-{\boldsymbol \nabla}_{i,t}, \x_{i,t}-x^\star_t} \leq 8LR \sqrt{-T\log \delta},
$$
with probability at least $1-\delta$. Summing \eqref{eq101} over $t\in [T]$ and $i\in [n]$ and incorporating the bound above into the last term, we get 
\begin{align*}
&\frac{1}{n}\sum_{i=1}^n\sum_{t=1}^Tf_{i,t}(\x_{i,t})-f_{i,t}(x^\star_t)\leq \\
&~~~~~~~~~~~~L^2\sum_{t=1}^T\frac{\eta_t}{2}+2L^2\sqrt{n}\sum_{t=1}^T\sum_{\tau=0}^{t-1}\eta_\tau\sigma^{t-\tau-1}_2(W)\\
\vphantom{\sum_{i=1}^n}&~~~~~~~~~~~~+\frac{1}{n}\sum_{i=1}^n\sum_{t=1}^T\left(\frac{1}{\eta_t}\dr(x^\star_t,\y_{i,t})-\frac{1}{\eta_t}\dr(x^\star_t,\xhx_{i,t+1})\right)\\
\vphantom{\sum_{i=1}^n}&~~~~~~~~~~~~+8LR \sqrt{-T\log \delta}, 
\end{align*}
with probability at least $1-\delta$. Applying Lemma \ref{auxlemma} to above, we can simplify as
\begin{align*}
&\frac{1}{n}\sum_{i=1}^n\sum_{t=1}^Tf_{i,t}(\x_{i,t})-f_{i,t}(x^\star_t)\leq \\
&~~~~~~~~~~~~L^2\sum_{t=1}^T\frac{\eta_t}{2}+2L^2\sqrt{n}\sum_{t=1}^T\sum_{\tau=0}^{t-1}\eta_\tau\sigma^{t-\tau-1}_2(W)\\
\vphantom{\sum_{i=1}^n}&~~~~~~~~~~~~+\frac{2R^2}{\eta_{T+1}}+\sum_{t=1}^T\frac{K}{\eta_{t+1}}\norm{x^\star_{t+1}-Ax^\star_{t}}\\
\vphantom{\sum_{i=1}^n}&~~~~~~~~~~~~+8LR \sqrt{-T\log \delta}. 
\end{align*}
We now return to sum \eqref{eq11} over $t\in [T]$. We apply the bound above and the bound in Lemma \ref{meandeviation}, respectively to the first and second line in \eqref{eq11} to finish the proof.
%\begin{align*}
%&\frac{1}{n}\sum_{i=1}^n\sum_{t=1}^Tf_t(\x_{i,t}) - f_t(x^\star_t)\leq \\
%&~~~~~~~~~~~~L^2\sum_{t=1}^T\frac{\eta_t}{2}+4L^2\sqrt{n}\sum_{t=1}^T\sum_{\tau=0}^{t-1}\eta_\tau\sigma^{t-\tau-1}_2(W)\\
%\vphantom{\sum_{i=1}^n}&~~~~~~~~~~~~+\frac{2R^2}{\eta_{T+1}}+\sum_{t=1}^T\frac{K}{\eta_{t+1}}\norm{x^\star_{t+1}-Ax^\star_{t}}\\
%\vphantom{\sum_{i=1}^n}&~~~~~~~~~~~~+8LR \sqrt{-T\log \delta}. 
%\end{align*}

\end{document}